\documentclass[11pt]{article}
\usepackage{amsmath, amsfonts, amssymb, amsthm}
\usepackage{natbib}

\DeclareMathSymbol{\leqslant}{\mathalpha}{AMSa}{"36} 

\DeclareMathSymbol{\geqslant}{\mathalpha}{AMSa}{"3E} 

\DeclareMathSymbol{\eset}{\mathalpha}{AMSb}{"3F}     

\newcommand{\suptwo}[2]{\sup_{\substack{#1 \\ #2}}} 

\newcommand{\supp}{\mathrm{\,range\,}}

\newcommand{\cC}{\mathcal{C}}

\newcommand{\cF}{\mathcal{F}}
\newcommand{\cG}{\mathcal{G}}
\newcommand{\cH}{\mathcal{H}}

\newcommand{\cL}{\mathcal{L}}
\newcommand{\cM}{\mathcal{M}}

\newcommand{\cO}{\mathcal{O}}

\renewcommand{\th}{\theta}
\newcommand{\Om}{\Omega}
\newcommand{\om}{\omega}

\newcommand{\ep}{\epsilon}
\newcommand{\eps}{\varepsilon}
\newcommand{\lam}{\lambda}
\newcommand{\sig}{\sigma}
\newcommand{\vphi}{\varphi}

\newcommand{\lfl}{\lfloor}
\newcommand{\rfl}{\rfloor}

\newcommand{\dto}{\downarrow}

\newcommand{\IP}{\mathbb{P}}
\newcommand{\II}{\mathbb{I}}
\newcommand{\IN}{\mathbb{N}}

\newcommand{\IC}{\mathbb{C}}
\newcommand{\IR}{\mathbb{R}}
\newcommand{\IE}{\mathbb{E}}
\newcommand{\iN}{\in\IN}

\newcommand{\iR}{\in\IR}

\newcommand{\be}{\begin{eqnarray*}}
\newcommand{\ee}{\end{eqnarray*}}
\newcommand{\ben}{\begin{eqnarray}}
\newcommand{\een}{\end{eqnarray}}

\theoremstyle{plain}
\newtheorem{theo}{Theorem}[section]
\newtheorem{lemma}[theo]{Lemma}
\newtheorem{propo}[theo]{Proposition}

\theoremstyle{definition}

\newtheorem{remark}[theo]{Remark}

\newtheorem{assu}[theo]{Assumption}

\renewenvironment{proof}[1][] {{\bf Proof#1.} }{\hspace*{\fill}$\square$\medskip\par}

\begin{document}
\vglue20pt
\centerline{\huge\bf The quantization complexity of}
\medskip

\centerline{\huge\bf  diffusion processes}

\bigskip
\bigskip

\centerline{by}
\bigskip
\medskip

\centerline{{\Large Steffen Dereich}}
\bigskip

\centerline{\it  Technische Universit\"at Berlin}

\bigskip
\bigskip
\bigskip

{\leftskip=1truecm
\rightskip=1truecm
\baselineskip=15pt
\small

\noindent{\slshape\bfseries Summary.} We investigate the high
resolution coding problem for solutions of stochastic differential equations in the $L^p[0,1]$-
and the $\IC[0,1]$-space. Tight asymptotic estimates are found under weak
regularity assumptions. The main technical tool is a decoupling method which allows
us to relate the complexity of the diffusion process to that of the Wiener
process in certain random Banach spaces.
\bigskip


\bigskip

}

\section{Introduction}\label{intro}

In this article, we study the high resolution quantization problem for $\IR^d$-valued stochastic processes $X$ that are solutions of stochastic differential equations. 
We let $\IC([0,1],\IR^d)$ denote the set of continuous functions  mapping $[0,1]$ to $\IR^d$ and denote by $\|\cdot\|_{[0,1]}$ the corresponding supremum norm, i.e.\ $\|f\|_{[0,1]}=\sup_{t\in[0,1]} |f(t)|$.
Analogously,  $L^q([0,1],\IR^d)$ denotes the $L^q$-space of functions mapping $[0,1]$ to $\IR^d$ and $\|\cdot\|_{L^q[0,1]}$ denotes the corresponding $L^q$-norm. 

We consider norm based distortion measures induced by
the latter two norms. 
More precisely, for $(E,\|\cdot\|)$ equal to $(\IC([0,1],\IR^d),\|\cdot\|_{[0,1]})$ or $(L^q([0,1],\IR^d), \|\cdot\|_{L^q[0,1]})$, we investigate the asymptotic behavior of the  quantization error
$$
D^{(q)}(r|X,\|\cdot\|,p)=\inf _{\hat X}\,  \IE[\|X-\hat X\|^p]^{1/p},
$$
where the infimum is taken over all discrete, $E$-valued random vectors $\hat X$ (reconstructions) with $|\supp(\hat X)|\le e^r$. 
The main objective is to provide sharp estimates for the asymptotic quantization error. Apart from the quantization problem, the new methods prove
to be useful also when considering entropy constrained coding.

Until the end of the past century, research was mainly
focused on source signals $X$  in finite dimensional spaces. An extensive
overview on this theory is contained in the monograph by
\cite{GraLu00}.  For a general absolutely continuous measure $\mu$ on $\IR^d$,
the asymptotic quantization error is related to that for the uniform
distribution on the unit cube; essentially, one obtains an additional
factor depending on the quantity $\|\frac{d\mu}{d\lam^d}\|_{L^q}$,
where $q>0$ is a parameter depending on the studied moment (see \cite{GraLu00}, Theorem 6.2).

In \cite{Feh01} and \cite{DFMS03}, the asymptotic quantization problem
was treated for Gaussian measures on separable (typically infinite dimensional) Banach spaces.
In this setting, the quantization error can be estimated against the
inverse of the small ball function. Results
from the theory of  small ball probabilities then lead to 
good estimates for many important examples. Moreover, in \cite{DerSch04} it is found that
for the $d$-dimensional Wiener process $W$ considered in $\IC([0,1],\IR^d)$, there exists a constant $K\in(0,\infty)$
independent of the moment index  $p>0$ such that the quantization error  satisfies
$$
D^{(q)}(r|W,\|\cdot\|_{[0,1]},p) \sim K\,\frac{1}{\sqrt r}
$$
as the rate $r$ tends to infinity. 
The constant $K$ is related to the principal eigenvalue of the Dirichlet problem on the unit disc. In particular, one obtains for $d=1$ that $K\in [\frac\pi{\sqrt {8}}, \pi ]$.

Here and elsewhere we write   $f \sim g$ iff
$\lim \frac fg = 1$,
while $f \lesssim g$ stands for $\limsup \frac fg \le 1$.
Finally, $f\approx g$ means
\[ 0< \liminf \frac fg \le \limsup \frac fg <\infty\ . \]


When the underlying space of the original process $X$ is a Hilbert space, 
a detailed analysis of the problem is possible: 
For any moment $p>0$ the asymptotic quantization error is
equivalent to the distortion-rate function for the second moment
norm-based distortion (mean squared error). The statement requires only
mild conditions on the
asymptotics of the eigenvalues of the corresponding covariance operator (see  \cite{Der03},
Theorem 6.2.1; \cite{LuPa04}). By a result of Kolmogorov (see \cite{Iha93}, Theorem 6.9.1) the distortion-rate function is given
by some implicit formula and the explicit asymptotics can be computed
in many cases. In particular, if the sequence of ordered eigenvalues $(\lam_n)$ is regularly varying with index $-\alpha<-1$, then for any moment index $p>0$,
$$
D^{(q)}(r|X, \|\cdot\|, p)\sim \Bigl( \frac{\alpha^\alpha}{2^{\alpha-1}(\alpha-1)}\, r \lam_{\lceil r\rceil}\Bigr)^{1/2}
$$
as $r\to\infty$ (see \cite{LuPa04}).
In particular, the $d$-dimensional Wiener process $W$ in $L^2([0,1],\IR^d)$
satisfies
$$
D^{(q)}(r|W,\|\cdot\|_{L^2[0,1]},2)\sim \frac{\sqrt2\, d}{\pi \sqrt r}.
$$

Let us now focus on the coding complexity of solutions $(X_t)_{t\in[0,1]}$ of stochastic differential equations. \cite{LuPa04b} considered $1$-dimensional diffusions with continuously differentiable diffusion coefficients. Their coding strategy is based on the Lamperti transform which maps the original $(X_t)$ onto a process $(\tilde X_t)$ which is a Brownian motion plus drift term. Approximating the process $\tilde X$ by some close process $\hat{\tilde X}$ and inverting the Lamperti transform for $\hat{\tilde X}$ leads to a "good" reconstruction of the original.
Under a regularity assumption on the Lamperti transform (assumption (3.8)) and the assumption that the diffusion coefficient is strictly bounded away from $0$, they are able to prove that
$$
D^{(q)}(r|X,\|\cdot\|_{L^q[0,1]},p)\approx \frac 1{\sqrt r},\qquad r\to\infty,
$$  
for any $q\in[1,\infty)$ and $p\in[1,\infty)$. The same estimate is valid when the $L^q$-norm is replaced by the supremum norm.

In contrast to \cite{LuPa04b}, we use the Doob-Meyer decomposition and a time change to approximate the original. Our approach leads to the strong asymptotics in the quantization problem when the underlying norm is the supremum norm. Moreover, we obtain an upper bound for the $L^q$-norm which we conjecture to be asymptotically tight. 
Our analysis requires only mild regularity assumptions for the drift and diffusion coefficients. Moreover, the multi-dimensional case with scalar diffusion coefficient is included in our approach. 
We shall see that the coding problem for  diffusions is tightly connected to that for the Wiener process, and this link will be the main tool in the proofs of the present results. Beyond the quantization problem, this fact seems to be useful as well when considering the entropy constrained coding problem.

Let us now fix the notation.
Let $(\Om,\cF,(\cF_t)_{t\ge0},\IP)$ be a complete filtered probability space that satisfies the usual hypotheses, i.e.\ $\cF_0$ contains all $\IP$-null sets of $\cF$ and $(\cF_t)$ is right continuous.
Let $(\tilde W_t)_{t\ge0}$ be a $d$-dimensional $(\cF_t)$-Wiener process. We denote by  $\sig:\IR^d\times[0,\infty)\to \IR$ and
$b:\IR^d\times[0,\infty)\to \IR^d$ two deterministic functions, 
and assume that $(X_t)_{t\ge0}$ is an $(\cF_t)$-adapted semimartingale solving the
integral equation
\begin{align}\label{eq0924-1}
X_t= \int_0^t b(X_s,s)\, ds +\int_0^t \sig(X_s,s)\ d\tilde W_s\qquad(t\ge0).
\end{align}
For ease of notation, we abridge $b_t:=b(X_t,t)$ and $\sig_t:=\sig(X_t,t)$ for $t\ge0$.
$(X_t)_{t\in[0,1]}$ represents the original process which will be approximated by some discrete r.v.\ $\hat X$, the reconstruction.

Our analysis requires the introduction of the quantization error under random distortion measures.
In the setting of a supremum norm-based distortions we  let $E=\IC([0,\infty),\IR^d)$, whereas $E=L^q([0,\infty),\IR^d)$ will be used for $L^q$-norm based distortions.
For an $E$-valued random vector $Z$, a measurable function $\rho:\Om\times
E\to[0,\infty)$, $p>0$, and $r\ge0$, let
\begin{align*}
D^{(q)}(r|Z,\rho, p)=\inf_{\hat Z} \,  \IE[\rho_{\om}(Z-\hat Z)^p]^{1/p},
\end{align*}
where the infimum is taken over all discrete, $E$-valued r.v.\ $\hat Z$ with $$|\supp (\hat
Z)|\le e^r.$$
This is the  \emph{$p$-th moment quantization error} for the
\emph{rate} $r$, \emph{source}  $Z$ and
\emph{distortion} $\rho$. 
An $E$-valued r.v.\ $\hat Z$ is called a \emph{regular reconstruction}
for $Z$, if
\begin{align}\label{eq2804-1}
\rho_\om(Z-\hat Z) \le \rho_\om(Z)
\end{align}
for all $\om\in\Om$.
Condition (\ref{eq2804-1}) is satisfied, if, for instance, $\hat Z$
is a reconstruction induced by a codebook containing the zero function $0$. The
quantization quantity obtained when confining oneself to regular reconstructions is denoted by $D^{(q,0)}(r|Z,\rho,p)$. 
Our analysis is based on a technical assumption:

\begin{assu}[C]\label{assuC}  There exist
  constants $\beta\in(0,1]$ and $L<\infty$ such that for
  $z,z'\iR^d\times [0,\infty)$
\begin{align}
|\sig(z)|+|b(z)| &\le L (|z|+1)\\
|\sig(z)-\sig(z')| & \le L [ |z-z'|^\beta + |z-z'|].
\end{align}
\end{assu}

Moreover, we assume that the  process $(\sig_t)_{t\in [0,1]}$ is not indistinguishible from the constant $0$-function.

Note that assumption (C) does neither ensure existence nor uniqueness of the solution of the stochastic differential equation (\ref{eq0924-1}).
A useful consequence of assumption (C) is that $\IE[\|X\|_{[0,1]}^p]$ is finite for any $p\ge1$.

Our main objectives are the following two theorems:

\begin{theo}\label{theo1109-1} Let  $p\ge1$ and 
$$
K=\lim_{r\to\infty} \sqrt r \, D^{(q)}(r|W,\|\cdot\|_{[0,1]},p)\in(0,\infty).
$$
Then 
$$
D^{(q)}(r|X,\|\cdot\|_{[0,1]},p)\sim K\,
\IE[\|\sig\|_{L^2[0,1]}^p]^{1/p}\,\frac1{\sqrt r}.
$$
\end{theo}

\begin{theo}
Let $p\ge1$ and $K<\infty$ be such that
$$
D^{(q)}(r|W,\|\cdot\|_{L^p[0,1]},p)\lesssim K\, \frac 1{\sqrt r}.
$$
Then 
$$
D^{(q)}(r|X,\|\cdot\|_{L^p[0,1]},p)\lesssim K\,
\IE[\|\sig\|_{L^{2p/(p+2)}[0,1]}^p]^{1/p}\,\frac1{\sqrt r}.
$$
\end{theo}

In order to approximate the process $(X_t)_{t\in[0,1]}$, we write $(X_t)$ in its Doob-Meyer decomposition $X_t=M_t+A_t$, where
$$
M_t=\int_0^t\sig(X_s,s)\,d\tilde W_s \ \text{ and }   \ A_t=\int_0^t b(X_s,s)\,ds.
$$
We shall see that the dominant term in the quantization problem is the continuous martingale $M$.
As is well known, we can represent  $M$ as a time change of a Wiener process. Let 
$$
\vphi(t)=\int_0^t \sig_s^2\,ds.
$$
By changing the drift and diffusion coefficients outside the time window $[0,1]$ we can ensure that $\lim_{t\to\infty}\vphi(t)=\infty$ without changing the process $(X_t)_{t\in[0,1]}$.
Then $W_t=
M_{\vphi^{-1}(t)}$  is a $d$-dimensional 
$(\cF_{t}^W)$-Wiener process, where
$$
\vphi^{-1}(t):=\inf\Bigl\{s\ge0: \int_0^s \sig^2_u\,du\ge t\Bigr\}.
$$

We roughly sketch the idea of the coding scheme for $M$. It can be decomposed in the following two steps:
\begin{enumerate}
\item approximate the real time transform $\vphi$ by some random monotone function $\hat\vphi\in\IC[0,1]$, and
\item approximate $(W_t)_{t\in[0,\tau]}$ ($\tau:= \hat\vphi(1)$) by $(\hat W_t)_{t\in[0,\tau]}$.
\end{enumerate}
Then $\hat M=\hat W_{\hat\vphi(\cdot)}$ is considered as the reconstruction, and the coding error can be controlled by
\begin{align*}
\|M-\hat M\|_{[0,1]} &\le
\|W_{\vphi(\cdot)}-W_{\hat\vphi(\cdot)}\|_{[0,1]}+\|W_{\hat
  \vphi(\cdot)}-\hat W_{\hat\vphi(\cdot)}\|_{[0,1]}\\
& \le \|W_{\vphi(\cdot)}-W_{\hat\vphi(\cdot)}\|_{[0,1]}+\|W
-\hat W\|_{[0,\tau]}.
\end{align*}

In the setting of $L^p$-norm based distortion,  
the corresponding estimate is then
\begin{align*}
\|M-\hat M\|_{L^p[0,1]}\le \|W_{\vphi(\cdot)}-W_{\hat\vphi(\cdot)}\|_{[0,1]}+\|W_{\hat\vphi(\cdot)}-\hat
W_{\hat\vphi(\cdot)}\|_{L^p[0,1]}.
\end{align*}
Denoting by $\nu$ the random measure induced by $\hat\vphi$, i.e.\ $\nu:=\lam(\hat\vphi^{-1})$, one
obtains
$$
\|M-\hat M\|_{L^p[0,1]}\le \|W_{\vphi(\cdot)}-W_{\hat\vphi(\cdot)}\|_{[0,1]}+\|W-\hat
W\|_{L^p(\nu)}.
$$
We shall see that, in both cases, the first term
in the sum is asymptotically negligible so that the asymptotics
are governed by the second terms. So we need strong estimates 
for the second term, whereas weak estimates suffice for the first term. 

The article is outlined as follows. The proofs of the theorems are based on a representation of the diffusion (Theorem
\ref{theo2705-1}) which will be given in Section \ref{sec7}.
The proof of this statement requires a couple of preliminary results:
We start in Section \ref{sec_entr} by providing an upper bound for the quantization
error based on entropy numbers of compact embeddings.
These yield weak asymptotic estimates which enable us to control
the asymptotically negligible terms.
Next, we prove an estimate for the
moments of the H\"older norm of continuous martingales. Based on this
estimate, we then show that the time transform $\vphi$ lies in a
H\"older-Zygmund space, and that all its moments are finite.
This leads to estimates for the
quantization problem for $\vphi$. The next section is devoted to an upper bound for
$\|W_{\vphi(\cdot)}-W_{\hat\vphi(\cdot)}\|_{[0,1]}$ for ``good''
reconstructions $\hat\vphi$ of $\vphi$.
In the following section, results of the theory of
enlargements of filtrations are used to show that roughly speaking the
martingale $W$ can be assumed to be independent of the time
transform $\hat\vphi$ without changing the coding complexity. Putting all these results together leads to the proof of the main representation (Theorem  \ref{theo2705-1}).
With this theorem, the problems of coding the diffusion are equivalent
to coding $W$ in  random $L^p(\nu)$ and $\IC[0,\tau]$-spaces, respectively. Solving
the optimal rate allocation problems leads to the main theorems.

\section{Entropy numbers and the quantization problem}\label{sec_entr}

Let $(E,\|\cdot\|_E)$ and $(F,\|\cdot\|_F)$ denote normed vector spaces such that $E$ is compactly
embedded into $F$. We denote by $e_n=e_n(E,F)$ the entropy numbers of
the embedding, i.e.
$$
e_n(E,F):=\inf\Bigl\{\eps>0: \exists x_1,\dots, x_{2^{n-1}}\in F \text{
  s.t.\ } B_E(0,1)\subset \bigcup_{i=1}^{2^{n-1}} B_F(x_i,\eps)\Bigr\}.
$$
We assume that $E$ is endowed with a $\sig$-field such that both norms
and the vector space operations are measurable.
In this section we write  $f\precsim g$ iff $\limsup \frac fg<\infty$.

\begin{lemma}\label{le0505-1} Let $\alpha>0$, and suppose that $E$ is compactly embedded into $F$
  with 
\begin{align}\label{eq0405-1}
e_n(E,F)\precsim n^{-\alpha},\qquad  n\to\infty .
\end{align}
Then for all $\tilde p>p>0$ there exists a constant
$c=c(p,\tilde p)<\infty$ such that  for all $E$-valued r.v.'s
$Z$ and $r\ge0$, one has
\begin{align}\label{eq20r}
D^{(q,0)}(r|Z,\|\cdot\|_F,p) \le c\, \IE[
\|Z\|_E^{\tilde p}]^{1/\tilde p}\,\frac1{1+r^\alpha}.
\end{align}
\end{lemma}

\begin{proof} Fix $\tilde p>p>0$. Notice that it suffices to prove
  the existence of a constant $c<\infty$  such that for any
  $E$-valued r.v.\ $Z$ with   $\IE[\|Z\|^{\tilde
    p}_E]^{1/\tilde p}= 1$:
$$
D^{(q,0)}(r|Z,\|\cdot\|_F,p) \le c\, \frac1{1+r^\alpha},
$$
since the general statement then follows by a scaling argument.

Notice that $e_n=e_n(E,F)$ is bounded by the norm
$\|\text{id}:E\to F\|=:\xi$. Consequently, using
assumption (\ref{eq0405-1}), there exists $c_1<\infty$ with
$$
e_n\le c_1 \,n^{-\alpha}\qquad(n\iN).
$$
Let $U=B_{E}(0,1)$ and 
$$
N(\eps,A)=\min\bigl\{|\cC|: A\subset \cC+ B_F(0,\eps)\bigr\}\qquad
(A\subset F, \, \eps>0).
$$
Then one has $\log N(2e_n,U)\le (n-1) \log 2$, and hence:  
$$
\log N(c_1 \,n^{-\alpha},U)\le (n-1) \log 2.
$$
By the monotonicity of $\eps\mapsto N(\eps,U)$, it now follows  that
$\log N(\eps,U)\precsim \eps^{-1/\alpha}$ as $\eps\dto0$. Since $N(\eps,U)=1$ for $\eps\ge
\xi$, there exists
a constant $c_2<\infty$ such that
$$
\log N(\eps,U)\le c_2\,\eps^{-1/\alpha}\qquad(\eps>0),
$$
and, consequently,
\begin{align}\label{eq22r}
\log N(\eps, sU)\le \frac{c_2\,s^{1/\alpha}}{\eps^{1/\alpha}}
\end{align}
for every $s,\eps>0$.
We fix $\eta>0$ such that $(1+\eta)p<\tilde p$, let $\eps>0$
arbitrary and consider
$$
\eps_i:=\eps_i(\eps):= \eps \,e^{(1+\eta) i}\qquad (i\iN_0).
$$
Moreover, let 
$$
s_i= e^i\qquad (i\iN_0),
$$
and $s_{-1}:=0$.
We use $ \eps_i$-nets of the sets $s_i U$ to generate an appropriate codebook.
Note that $\eps_i\ge \xi s_i$, if 
\begin{align*}
i \ge \Bigl\lceil \frac 1\eta \log(\xi/\eps)\Bigr\rceil \vee 0 =:M.
\end{align*}
Since $\xi \|x\|_E \ge \|x\|_F$ for $x\in E$, the set $\{0\}$ is an optimal
$\eps_i$-net of $r_i U$ for $i\ge M$.
We consider the codebook
$$
\cC(\eps)= \{0\}\cup \bigcup_{i=0}^{M-1} \cC_i(\eps),
$$
where $\cC_i(\eps)$ denote arbitrary optimal $\eps_i$-nets of $r_i U$ $(i\iN_0)$. Then
\begin{align*}
\IE[d_F(Z,\cC(\eps))^p] & \le  
\sum_{i=0}^\infty \IE[1_{[s_{i-1}, s_{i})}(\|Z\|_{E})  \,d(Z,\cC_i(\eps))^p]\\
& \le \sum_{i=0}^\infty \IP(\|Z\|_{E}\ge s_{i-1})\, \eps_i^p\\
& = \eps^p +\sum_{i=1}^\infty \IP\Bigl(\frac{\|Z\|_{E}^{\tilde p}}{s_{i-1}^{\tilde
    p}}\ge 1\Bigr)\, \eps_i^p\\
& \le \eps^p + \IE[\|Z\|_{E}^{\tilde p}]\,\sum_{i=1}^\infty \frac{\eps_i^p}{s_{i-1}^{\tilde
    p}}\\
& = \eps^p \Bigl(1+
\,\sum_{i=1}^\infty \frac{e^{(1+\eta)ip}}{e^{(i-1)\tilde
    p}}\Bigr)\\
&= \eps^p \Bigl(1+
\,\sum_{i=1}^\infty e^{\tilde p-(\tilde p-(1+\eta)p) i}\Bigr).
\end{align*}
Since $\tilde p>(1+\eta)p$, the previous sum converges.
Consequently, there exists a constant $c_3<\infty$ not depending on 
$Z$ and such that
$$
\IE[d_F(Z,\cC(\eps))^p]  \le c_3\,\eps^p.
$$
It remains to compute an upper bound for the size of $\cC(\eps)$. If $\eps\ge \xi$, then $M= 0$ and
$|\cC(\eps)|=1$. On the other hand, for $\eps<\xi$  equation
(\ref{eq22r}) implies
\begin{align*}
|\cC(\eps)| &\le 1+ \sum_{i=0}^{M-1} |\cC_i(\eps)| \le 1+
\sum_{i=0}^{M-1} \exp\bigl\{ c_2 \,
(s_i/\eps_i)^{1/\alpha}\bigr\}\\
& = 1+
\sum_{i=0}^{M-1} \exp\bigl\{ c_2 \,\frac1{\eps^{1/\alpha}}
\,e^{-i\eta/\alpha}\bigr\}\\
& \le  1+  M\, \exp\bigl\{ c_2
\,\frac1{\eps^{1/\alpha}}\bigr\}\\
& \le 1+ \Bigl(1+\frac1\eta\,\log (\xi/\eps)\Bigr)\, \exp\bigl\{ c_2
\,\frac1{\eps^{1/\alpha}}\bigr\} .
\end{align*}
Combining both estimates, we conclude that there exists a constant $c_4<\infty$ such that for all $\eps>0$,
$$
|\cC(\eps)|\le \exp\bigl\{c_4 \,\frac1{\eps^{1/\alpha}}\bigr\}.
$$
Therefore,
$$
D^{(q,0)}\bigl(c_4\,\eps^{-1/\alpha}\,|Z,\|\cdot\|_F,p\bigr)^p \le c_3\,\eps^p
$$
for $\eps>0$ and, hence, choosing $r>0$ arbitrary and $\eps=(r/c_4)^{-1/\alpha}$ yields 
\begin{align*}
D^{(q,0)}(r|Z,\|\cdot\|_F,p)\le c_5\,\frac1{r^\alpha}
\end{align*}
for some constant $c_5<\infty$. Since $D^{(q,0)}(r,|Z,\|\cdot\|_F,p)$
is uniformly bounded by $\xi$, we finally conclude that there exists a
constant $c_6<\infty$ such that
\begin{align*}
D^{(q,0)}(r|Z,\|\cdot\|_F,p)\le c_6\,\frac1{1+r^\alpha}
\end{align*}
for all $r\ge0$.
\end{proof}

\section{H\"older continuity of $M$}

Let  $M=(M_t)_{t\in[0,1]}$ be an $(\cF_t)$-adapted process of the form $M_t=\int_0^t
\sigma_s\,d\tilde W_s$, where $(\sigma_t)$ is an  $(\cF_t)$-adapted
process such that the integral is well defined. In this
section we do not require that $(\sigma_t)$ be given by $\sigma_t=\sigma(X_t,t)$.   

We denote by $|\cdot|_\alpha$ the $\alpha$-H\"older semi-norm on
$\IC[0,1]$, i.e.
$$
|f|_\alpha:=\sup_{0\le s<t\le1} \frac {|f(t)-f(s)|}{|t-s|^{\alpha}}.
$$

Based on the GRR inequality (see \cite{GRR71}) we derive an upper bound
for the moments of $|M|_\alpha$:

\begin{theo}\label{theo0806-1}
Let $\alpha\in(0,1/2)$ and $\kappa>2/(1-2\alpha)$. Then there exists a
constant $c=c(\kappa,\alpha)$ such that
$$
\IE[|M|^\kappa_\alpha]\le c\,\int_0^1 \IE[|\sig_u|^\kappa]\,du.
$$
The constant $c$ does not depend on the martingale $M$.
\end{theo}

\begin{proof}
Fix $\alpha\in(0,1/2)$.
Let
$f:[0,1]\to \IR$ be a continuous function and let $\beta,\gamma>0$ with
$\alpha=\gamma-2/\beta$. We consider  $\Psi(x)=|x|^\beta$ and
$p(x)=|x|^\gamma$, $x\iR$.
Then the GRR lemma states that for
\begin{align}\label{eq2101-2}
B:=B(f):=\int_0^1\int_0^1\Psi\Bigl( \frac {f(s)-f(t)}{p(s-t)}\Bigr)\,ds\,dt=\int_0^1\int_0^1\frac {|f(s)-f(t)|^\beta }{|s-t|^{\beta \gamma}}\,ds\,dt,
\end{align}
one has
\begin{align*}
|f(s)-f(t)|&\le 8 \int_0^{|s-t|}\Psi^{-1}\Bigl(\frac{4B}{\xi^{2}}\Bigr)\,dp(\xi)\\
&=  8
\int_0^{|s-t|}\frac{(4B)^{1/\beta}}{\xi^{2/\beta}}\,dp(\xi)\\
&=8\gamma (4B)^{1/\beta}\int_0^{|s-t|}\xi^{\gamma-1-2/\beta} \,d\xi\\
&=8\frac{\gamma}{\gamma-2/\beta} (4B)^{1/\beta} |s-t|^{\gamma-2/\beta}
\end{align*}
for all $s,t\in[0,1]$.
Consequently,
$$
|f|_\alpha \le 4^{1/\beta} 8\frac{\gamma}{\alpha} B^{1/\beta}.
$$

Now fix $i\in\{1,\dots,d\}$ and consider the $i$-th coordinate of $M$
(denoted by $M^i$) instead of $f$. We define $B:=B(M^i)$ in analogy to
(\ref{eq2101-2}), and estimate the $\kappa$-th moment of $|M^i|_\alpha$ for $\kappa\ge\beta\vee2$.
Using Jensen's inequality and the  Burkholder-Davis-Gundy (BDG)
inequality, we conclude that there exists some constant
$c_1=c_1(\kappa)<\infty$ such that
\begin{align*}
\IE[ B^{\kappa/\beta}]&\le \IE\Bigl[ \int_0^1\int_0^1\frac
{|M^i_s-M^i_t|^\kappa }{|s-t|^{\kappa \gamma}}\,ds\,dt\Bigr]\\
&\le  \int_0^1\int_0^1 \frac
{\IE[|M_s-M_t|^\kappa] }{|s-t|^{\kappa \gamma}}]\,ds\,dt\\
&\le \int_0^1\int_0^1 \frac
{c_1 \,\IE[|\int_s^t |\sig_u|^2 du|^{\kappa/2}]
}{|s-t|^{\kappa \gamma}}\,ds\,dt\\
&= 2 c_1 \int_0^1\int_0^{1-\delta} \frac
{\IE[(\int_s^{s+\delta} |\sig_u|^2 du)^{\kappa/2}]
}{\delta^{\kappa \gamma}}\,ds\,d\delta.
\end{align*}
Applying again the Jensen inequality leads to
\begin{align*}
\IE[ B^{\kappa/\beta}]&\le 2 c_1 \int_0^1\delta^{-\kappa\gamma +\kappa/2}\int_0^{1-\delta}
 \IE\Bigl[ \int_s^{s+\delta} \delta^{-1} |\sig_u|^\kappa du \Bigr]
 \,ds\,d\delta.
\end{align*}
By elementary analysis, one can show that the inner two integrals
satisfy, for $\delta\in(0,1)$:
\begin{align*}
\int_0^{1-\delta}
 \int_s^{s+\delta} \delta^{-1} \IE[|\sig_u|^\kappa]\, du
 \,ds\le \int_0^1\IE[|\sig_u|^\kappa]\,du.
\end{align*}
Consequently,
$$
\IE[ B^{\kappa/\beta}]\le 2 c_1 \int_0^1\delta^{-\kappa\gamma
  +\kappa/2}\,d\delta \cdot  \int_0^1\IE[|\sig_u|^\kappa]\,du.
$$
Thus if $-\kappa\gamma+\kappa/2>-1$, then there exists a constant
$c_2<\infty$ depending on $\kappa$ and $\gamma$ with
\begin{align}\label{eq2101-1}
\IE[|M^i|_\alpha ^\kappa]\le c_2\,\int_0^1\IE[|\sig_u|^\kappa]\,du.
\end{align}
Now we need to study for which $\kappa\iR_+$ there exist appropriate
values for $\beta$ and $\gamma$ admitting an estimate (\ref{eq2101-1})
with finite $c_2$: $\beta$ and $\gamma$ need to satisfy
$$
(i)\  \gamma-2/\beta=\alpha\ ,\ (ii)\ \kappa\ge \beta\vee 2\ \text{ and }\
(iii)\ -\kappa \gamma+\kappa/2>-1.
$$
When choosing $\gamma\in (1/2,\infty)$, conditions (i) to (iii) are
equivalent to 
$$
\kappa\ge \frac2{\gamma-\alpha}\vee 2\ \text{ and } \ \kappa<\frac1{\gamma-1/2}.
$$
Elementary analysis implies the existence of an estimate like (\ref{eq2101-1})  for each 
$$
\kappa>\frac2{1-2\alpha}.
$$
The general assertion of the lemma is obtained via applications of the
triangle inequality.
\end{proof}

\begin{remark}
The condition that $\kappa>2/(1-2\alpha)$ is necessary for the
validity of the lemma. If the condition is not satisfied, a counterexample is
obtained as follows: fix  $\eps\in(0,1]$ and let $\sig_t:=1_{[0,\eps]}(t)
\eps^{-1/\kappa}$ ($t\in[0,1]$); then the right hand side of the inequality is equal
to $c$, whereas  $\IE[|M|_\alpha^\kappa]$ tends to infinity
when letting $\eps$ to zero.
\end{remark}

\section{Regularity of  $\vphi$ and its coding complexity}

We assume again the setting of section \ref{intro}. 
Let $m\iN$, $\alpha\in(0,1]$ and $s:=m+\alpha$. For 
$f\in \IC^m[0,1]$ we let
$$
\|f\|_s:=\|f\|_{[0,1]}+|f^{(m)}|_\alpha
$$
denote the $s$-th H\"older-Zygmund norm and denote by 
$$
C^s:=\{f\in \IC^m[0,1]: \|f\|_s<\infty\}
$$
the corresponding Banach space.
Owing to Kolmogorov, the space $C^s$ is compactly embedded into
$\IC[0,1]$ and the corresponding  metric entropy numbers satisfy
$$
e_n(C^s,\IC[0,1])\approx \frac 1{n^s},\qquad n\to\infty.
$$
Thus Lemma \ref{le0505-1} immediately implies:

\begin{lemma} \label{le0906-1} Let $\tilde p>p>0$ and $s>0$. There exists a constant
  $c=c(p,\tilde p,s)$ such that for all  $C^s$-valued
  random variables $Z$ and all $r\ge0$
$$
D^{(q,0)}(r|Z,\|\cdot\|_{[0,1]},p)\le c \,\IE[\|Z\|_s^{\tilde p}]^{1/\tilde p}\,\frac1{r^s+1}.
$$
\end{lemma}

We use this fact to prove the following lemma.

\begin{lemma}  \label{le1604-2} For $p>0$ and $\alpha\in(0,\beta/2)$,
  there exists a constant   $c<\infty$ such that
$$
D^{(q,0)}(r|\vphi,\|\cdot\|_{[0,1]},p)\le \frac c{r^{1+\alpha}+1}.
$$
\end{lemma}

\begin{proof}
Recall that $\vphi(t)=\int_0^t \sig_s^2\,ds$. Thus, in view of Lemma \ref{le0906-1}, it suffices to prove
  that for $p\ge1$ and $\alpha\in(0,\beta/2)$,
$$
\IE[\|\vphi\|_{[0,1]}^p]^{1/p}+ \IE[|\sig^2|_\alpha^p]^{1/p}<\infty.
$$
Note that by assumption (C):
$$
\IE[\|\sig\|_{[0,1]}^p]^{1/p}\le 2L +L\, \IE[\|X\|_{[0,1]}^p]^{1/p} <\infty.
$$
In particular, $\IE[\|\vphi\|_{[0,1]}^p]$ is finite for all $p\ge1$. It remains to consider $\IE[|\sig^2|_{\alpha}^p]$.
By Theorem \ref{theo0806-1}, it is true that for every $\alpha'\in(0,1/2)$ and $p\ge1$
$$
\IE[|M|_{\alpha'}^p]<\infty.
$$
Moreover, again by assumption (C)
$$
\IE[|A|_{\alpha'}^p]^{1/p}\le \IE[|A|_1^p]^{1/p}\le \IE[\|b\|_{[0,1]}^p]^{1/p}\le 2L+L \ \IE[\|X\|_{[0,1]}^p]^{1/p}<\infty
$$
and hence 
\begin{align}\label{eq1027-1}
\IE[|X|_{\alpha'}^p]^{1/p} \le \IE[|M|_{\alpha'}^p]^{1/p} +\IE[|A|_{\alpha'}^p]^{1/p}<\infty. 
\end{align}
Since in general $|\sig^2|_{\alpha}\le 2\|\sig\|_{[0,1]} |\sig|_\alpha$,
one also has:
$$
\IE[(|\sig^2|_{\alpha})^p]\le 2^p\ \IE[\|\sig\|_{[0,1]}^p |\sig|_{\alpha}^p],
$$
and due to the Cauchy Schwarz inequality it suffices to establish the finiteness of $\IE[|\sig|_{\alpha}^{2p}]$.
By elementary analysis, we obtain that
$$
|\sig|_\alpha\le L \bigl( |X|_{\alpha/\beta}^{\beta }+ |X|_{\alpha}+2\bigr),
$$
so that  (\ref{eq1027-1}) and the inequality $\alpha/\beta<1/2$ imply that $\IE [|\sig|_\alpha^{2p}]$ is finite. 
\end{proof}

\section{An estimate for $\IE[\|W_{\vphi(\cdot)}-W_{\hat\vphi(\cdot)}\|_{[0,1]}^p]$}

In the previous section we considered bounds for the quantization
problem for the original $\vphi$ in $\IC[0,1]$. It remains to study
the quantity
$\IE[\|W_{\vphi(\cdot)}-W_{\hat\vphi(\cdot)}\|_{[0,1]}^p]$ for
``good'' reconstructions $\hat \vphi$ of $\vphi$. The following
analysis relies severely on concentration properties of Gaussian measures.

\begin{lemma}\label{le2}
Let $T$, $\eps_1$, $\eps_2>0$ with $\eps_2\ge\sqrt{2\eps_1}$. We have
$$
\IP\Bigl(\suptwo{\{s,t\in[0,T]:}{|s-t|\le \eps_1\}} |W_t-W_s|\le 3\eps_2\Bigr)\ge\Bigl(1-2\,e^{-\frac{\eps_2^2}{2\eps_1}}\Bigr)^n,
$$
where $n:=\lceil T/\eps_1\rceil$.
\end{lemma}

\begin{proof}
Let $T, \eps_1$, $\eps_2$ and $n$ as in the lemma.
Set $t_i=i\eps_1$, $i=0,\dots, n-1$, and $t_n=T$.
It holds
$$
\suptwo{\{s,t\in[0,T]:}{|s-t|\le \eps_1\}} |W_t-W_s|\le 3
\max_{i=0,\dots, n-1} \sup_{s\in [t_i,t_{i+1}]} |W_s-W_{t_i}|.
$$
We denote $M_i=\sup_{s\in [t_i,t_{i+1}]} |W_s-W_{t_i}|$,
$i=0,\dots,n-1$.
Then
$$
\IP\Bigl(\suptwo{\{s,t\in[0,T]:}{|s-t|\le \eps_1\}} |W_t-W_s|\le 3\eps_2\Bigr)\ge \IP(\max_{i=0,\dots, n-1} M_i\le \eps_2)
$$
Note that the random variables $M_0, \dots, M_{n-1}$ are
independent and, therefore,
$$
\IP(\max_{i=0,\dots, n-1} M_i\le \eps_2)=\prod_{i=0}^{n-1} \IP( M_i\le \eps_2).
$$
We compute
\begin{align*}
\IP( M_i> \eps_2)& \le 2\,\IP(\sup_{s\in [t_i,t_{i+1}]}
(W_s-W_{t_i})>\eps_2)\\
& = 4 \,\IP(W_{t_{i+1}}-W_{t_i})>\eps_2)\\
& = 4 \,\IP\Bigl(\frac1{\sqrt{t_{i+1}-t_i}} (W_{t_{i+1}}-W_{t_i})>\frac {\eps_2}{\sqrt{t_{i+1}-t_i}}\Bigr)\\
& = 4\, \Psi\Bigl(\frac {\eps_2}{\sqrt{t_{i+1}-t_i}}\Bigr)\le 2\,e^{-\frac{\eps_2^2}{2(t_{i+1}-t_i)}},
\end{align*}
where $\Psi(t):=(2\pi)^{-1/2}\int_t^\infty \exp\{-x^2/2\}\,dx$ ($t\iR$).
By assumption, the last term is less than $1$ and  
$$
\IP(\max_{i=0,\dots, n-1} M_i\le \eps_2)\ge \prod_{i=1}^{n-1}
\Bigl(1-2\,e^{-\frac{\eps_2^2}{2(t_{i+1}-t_i)}}\Bigr)\ge \Bigl(1-2\,e^{-\frac{\eps_2^2}{2\eps_1}}\Bigr)^n.
$$
\end{proof}

Notice that $\IC_0[0,T]$ equipped with the norm
$$
\|f\|_{\eps_1,T}:=\suptwo{s,t\in[0,T]:}{|s-t|\le \eps_1} |f(t)-f(s)|
$$
is a separable Banach space, say $C_{\eps_1,T}$. Thus we can interpret $(W_t)_{t\in[0,T]}$ as a
centered Gaussian random vector in
this space.  Let $m_{\eps_1,T}\iR_+$ denote the $7/8$-quantile of
$\|W\|_{\eps_1,T}$. 
Using  elementary analysis together with Lemma \ref{le2}, we obtain:

\begin{lemma} \label{le0303-2}
There exists  a constant $c<\infty$ such that for
  all $T,\eps_1>0$ one has
$$
m_{\eps_1,T}\le \tilde c\, \sqrt{\eps_1 \Bigl(1+\log_+ \frac T{\eps_1}\Bigr)}.
$$ 
\end{lemma}

\begin{lemma}
For any $A\in\cF$ and $p>0$, it holds
$$
\IE[1_A \ \|W\|_{\eps_1,T}^p]^{1/p}\le 2\sqrt2\ m_{\eps_1,T}
\Bigl[2\int^\infty_{0\vee \Psi^{-1}(\IP(A))}
\frac{(x+1)^p}{\sqrt{2\pi}}\,e^{-x^2/2}\ dx\Bigr]^{1/p},
$$
where $\Psi(t):=\int_t^\infty \sqrt{2\pi}^{-1} \exp\{-x^2/2\}\,dx$ $(t\iR)$.
\end{lemma}

\begin{proof}
Due to \cite{LedTal91} (p.\ 99) (see also \cite{Led96}, p.\ 202, 210)  it holds
$$
\sig:=\sup_{f\in C_{\eps_1,T}^*, \|f\|_{C_{\eps_1,T}^*}\le1} \IE[f^2(W)]^{1/2}\le 2\sqrt2\ m_{\eps_1,T},
$$
where $C^*_{\eps_1,T}$ is the topological dual of $C_{\eps_1,T}$.
As a consequence of the isoperimetric inequality, one obtains
$$
\IP(\|W\|_{\eps_1,T}\ge m_{\eps_1,T}+ t \sig)\le \Psi(t)\qquad (t\ge0).
$$
Therefore, we can find a standard normal
random variable $N$ on a sufficiently large probability space such that
$$
\|W\|_{\eps_1,T}\le 2\sqrt2\ m_{\eps_1,T}[1+N^+],
$$
where $N^+=N\vee 0$.
Consequently,
$$
\IE[1_A \ \|W\|_{\eps_1,T}^p]^{1/p}\le 2\sqrt2\
m_{\eps_1,T}\,\IE[1_A\, (N^++1)^p]^{1/p}.
$$
\end{proof}

\begin{lemma} \label{le0504-1} For $p>0$ there exists a constant $c=c(p)$
  such that for all $A\in\cF$, $T,\eps_1>0$ one has
$$
\IE[1_A\|W\|^p_{\eps_1,T}]^{1/p}\le c\,\sqrt{\eps_1 \Bigl(1+\log_+ \frac T{\eps_1}\Bigr)}\,q\,\sqrt{1+\log(1/q)}^{p-1},
$$
where $q:=\IP(A)$. 
\end{lemma}

\begin{proof}By elementary analysis one obtains
$$
\int_x^\infty (y+1)^p\, e^{-y^2/2}\,dy\sim x^{p-1}\,\Psi(x),\qquad x\to\infty,
$$
and thus
$$
\int_{\Psi^{-1}(\eps)}^\infty (y+1)^p\, e^{-y^2/2}\,dy\sim \eps\,\Psi^{-1}(\eps)^{p-1}\sim \eps\,\sqrt{2\log(1/\eps)}^{p-1},\qquad \eps\dto0.
$$
Consequently, there exists some constant $c_1=c_1(p)<\infty$ such
that for all $\eps\in(0,1]$
$$
\int_{0\vee \Psi^{-1}(\eps)}^\infty (y+1)^p\, e^{-y^2/2}\,dy\le
c_1 \,\eps\,\sqrt{1+\log(1/\eps)}^{p-1}
$$
An application of the previous two lemmas yields
$$
\IE[1_A\|W\|^p_{\eps_1,T}]^{1/p}\le \kappa_p''\,\sqrt{\eps_1 \Bigl(1+\log_+ \frac T{\eps_1}\Bigr)}\,q\,\sqrt{1+\log(1/q)}^{p-1},
$$
where $q:=\IP(A)$ and $c_2=c_2(p)$ is a constant only depending on $p$.
\end{proof}

\begin{lemma} \label{le1604-3} Suppose  that $\hat\vphi^{(r)}$ ($r\ge0$) are reconstructions for
$\vphi$ such that
$$
\lim_{r\to\infty} \IE[\|\vphi-\hat \vphi^{(r)}\|_{[0,1]}^2]^{1/2}=0.
$$
Then for any $p\ge1$,
$$
\IE[\|W_{\vphi(\cdot)}-W_{\hat\vphi^{(r)}(\cdot)}\|_{[0,1]}^p]^{1/p}=
\cO\left(\sqrt{d(r) \,\log(1/d(r))}\right),
$$
where $d(r)=\IE[\|\vphi-\hat \vphi^{(r)}\|_{[0,1]}^2]^{1/2}$.
\end{lemma}

\begin{proof} 
Consider the r.v.'s
  $\eps_1:=\eps_1(r):=\|\vphi-\hat\vphi^{(r)}\|_{[0,1]}$ and
  $\tau:=\sup_{t\in[0,1]}\vphi(t)$.
Notice that
$$
\|W_{\vphi(\cdot)}-W_{\hat\vphi^{(r)}(\cdot)}\|_{[0,1]}\le \|W\|_{\eps_1,\tau+\eps_1}.
$$
Let now $\II:=\{e^i:i\iN_0\}$,
$$
\bar \eps_1:=\bar \eps_1(r):=\min( [\eps_1,\infty)\cap \IE[\eps_1^2]^{1/2}\II) \text{
  and }\bar \tau :=\bar \tau (r):=\min( [\eps_1,\infty)\cap \II).
$$
$\bar\eps_1$ and $\bar \tau$ are discrete r.v.'s dominating
$\eps_1$ and $\tau$ and satisfying
\begin{align}\label{eq2804-2}
\bar\eps_1\le e\eps_1+\IE[\eps_1^2]^{1/2}\text{ and }\bar\tau\le e\tau+1.
\end{align}
Denote by
$(p_{\eps,t})$ the probability weights of $(\bar \eps_1,\bar
\tau)$. Then the triangle inequality and Lemma \ref{le0504-1} yield
\begin{align*}
\IE[\|W\|_{\eps_1,\tau+\eps_1}^p]^{1/p}&\le
\IE[\|W\|_{\bar\eps_1,\bar\tau+\bar\eps_1}^p]^{1/p}\\
&\le \sum_{\eps,t} \IE[ 1_{\{(\eps_1,\tau)=(\eps,t)\}}
\|W\|_{\eps,t+\eps}^p]^{1/p}\\
&\le c_1 \,\IE\Bigl[ \sqrt{\bar\eps_1(1+\log(1+\frac {\bar\tau}{\bar\eps_1}))}
 \sqrt{1+\log(1/p_{\bar\eps_1,\bar\tau})}^{p-1}\Bigr]\\
&\le c_1 \, \IE\Bigl[ \bar\eps_1(1+\log(1+\frac
{\bar\tau}{\bar\eps_1}))\Bigr]^{1/2}
\,\IE\Bigl[\bigl(1+\log(1/p_{\bar\eps_1,\bar\tau})\bigr)^{p-1}\Bigr]^{1/2}\\
&=:c_1\, \Sigma_1 \cdot \Sigma_2
\end{align*}
for some appropriate constant $c_1=c_1(p)$.
Notice that the second term is dominated by 
$$
\Sigma_2\le c_2\, (H^{p-1}(\bar\eps_1,\bar\tau)^{1/2}+1),
$$
where $c_2=c_2(p)<\infty$ is a universal constant and $H^{p-1}$ denotes the
generalized entropy 
$$
H^{p-1}(\bar\eps_1,\bar\tau):=\sum_{\eps,t} p_{\eps,t}\bigl(\log(1/p_{\eps,t})\bigr)^{p-1}.
$$
Now choose $\eps=\IE[\eps_1^2]^{1/2}e^i$ and $t=e^j$ ($i,j\iN_0$). If
$i,j\iN$, one obtains with (\ref{eq2804-2}) and the Cauchy-Schwarz
inequality 
\begin{align*}
p_{\eps,t}&\le \frac {\IE[\bar\eps_1
  \bar\tau]}{\IE[\eps_1^2]^{1/2}e^{i+j-2}}\le \frac
{\IE[(e\eps_1+\IE[\eps_1^2]^{1/2})(\tau+1)]}{\IE[\eps_1^2]^{1/2}e^{i+j-2}}\\
&\le \frac{(e+1)\, \IE[(\tau+1)^2]^{1/2}}{e^{i+j-2}}.
\end{align*}
If $i=0$ and $j\iN$, then
\begin{align*}
p_{\eps,t}&\le \frac {\IE[ \bar\tau]}{e^{j-1}}\le \frac
{\IE[(\tau+1)]}{e^{j-1}}\le \frac{ \IE[(\tau+1)^2]^{1/2}}{e^{j-1}},
\end{align*}
whereas for $i\iN$ and $j=0$, one obtains
$$
p_{\eps,t}\le \frac {\IE[ \bar\eps_1]}{\IE[\eps_1^2]^{1/2} e^{i-1}}
\le \frac1{e^{i-1}}.
$$
Note that the above estimates for $p_{\eps,t}$ do not depend on the rate $r\ge0$ and decrease sufficiently fast to zero in order to provide the finiteness of $H^{p-1}(\bar\eps_1,\bar\tau)$. Therefore, $H^{p-1}(\bar\eps_1,\bar\tau)$ is uniformly bounded for all $r\ge0$ by some constant $c_3<\infty$ depending only on $\IE[(\tau+1)^2]^{1/2}$. Consequently, $\Sigma_2$ is uniformly bounded.

It remains to consider the first expression $\Sigma_1$. Using (\ref{eq2804-2})
and the inequality $\bar\tau\ge1$ we arrive at
\begin{align*}
\Sigma_1^2&= \IE\Bigl[ \bar\eps_1(1+\log(1+\frac
{\bar\tau}{\bar\eps_1}))\Bigr]\le \IE\Bigl[ \bar\eps_1(1+\log\bar \tau+
\log(1+\frac
{1}{\bar\eps_1}))\Bigr]\\
&\le\IE\Bigl[(e\eps_1+\IE[\eps_1^2]^{1/2})\Bigl(1+\log(1+e\tau) +\log \Bigl(1+\frac{1}{e\eps_1+\IE[\eps_1^2]^{1/2}}\Bigr)\Bigr)\Bigr].
\end{align*}
An application of the Cauchy-Schwarz inequality then yields that
$$
\Sigma_1^2=\cO\bigl(\IE[\eps_1^2]^{1/2} \log(1/\IE[\eps_1^2]^{1/2})\bigr)
$$
and the assertion follows.
\end{proof}

\section{Coding $(W_t)_{t\in[0,\tau]}$}

We need some more notations. For $f\in \IC_0([0,\infty),\IR^d)$, let
\begin{align*}
\|f\|_{\cH}:=\begin{cases}
\|\frac {df}{dt}\|_{L^2[0,\infty)}  & \text{if }f\text{ is
  weakly differentiable}\\
\infty & \text{else}.
\end{cases}
\end{align*}
Moreover, for $T>0$ and $f\in \IC_0([0,T],\IR^d)$, let $\|f\|_{\cH_T}=\|\frac {df}{dt}\|_{L^2[0,T]}$ if $f$ is weakly
differentiable on $[0,T]$ and $\|f\|_{\cH_T}=\infty$, otherwise. The corresponding
Hilbert spaces are denoted by $\cH$ and $\cH_T$.

We recall some results of the theory of enlargements of filtrations
(see \cite{JeuYor85}, \cite{AnkDerImk05}).
Let $(\cF^W_t)$ be the filtration generated by the Wiener process
$(W_t)$ and denote by $Z$ a discrete  random variable with probability
weights $(p_x)$.
We consider the enlarged filtration $\cG_t=\cF^W_t\vee \sig(Z)$
$(t\ge0)$ and assume that for some fixed $p\ge1$ the generalized entropy
$$
H^p(Z):=\IE\Bigl[\Bigl(\log \frac 1{p_Z}\Bigr)^p\Bigr]
$$
is finite.
Then the process $(W_t)$ is a $(\cG_t)$-semimartingale, and
its Doob-Meyer decomposition  $W_t=\bar W_t+\bar A_t$  comprises a
$(\cG_t)$-Wiener process $(\bar W_t)$ and a process of bounded variation
$(\bar A_t)$ satisfying
\begin{align}\label{eq2204-1}
\IE[\|\bar A\|^{2p}_\cH]\le
\kappa_p \ (H^{p}(Z)+1).
\end{align}
Here, the constant $\kappa_p$ depends only on $p$.

We recall that  $\cH_1$ is compactly embedded into $\IC([0,1],\IR^d)$ and that its
entropy numbers satisfy
$$
e_n(\cH_1,\IC([0,1],\IR^d))\approx \frac1n,\qquad n\to\infty.
$$

\begin{lemma}\label{le0705-1}
Let  $\tilde p>p>0$. There exists a constant
  $c=c(p,\tilde p)$ such that 
$$
D^{(q,0)}(r|Y,\|\cdot\|_{[0,T]},p) \le c\ \sqrt T\
\IE[\|Y\|_{\cH_T}^{\tilde p}]^{1/\tilde p}\ \frac 1{r+1}
$$ 
for all $T>0$, $r\ge0$ and $\cH_T$-valued r.v.\ $Y$. 
\end{lemma}

\begin{proof}
By Lemma \ref{le0505-1}, the statement holds for fixed time $T=1$ for an appropriate constant $c>0$. 
Notice that for $T>0$ the maps
\begin{align*}
\pi^{(1)}_T&:\cH_T\to \cH_1,\, f\mapsto \frac1{\sqrt T} f(T\cdot)\quad\text{
  and}\\
 \pi^{(2)}_T&:\IC[0,T]\to \IC[0,1],\, f\mapsto  f(T\cdot)
\end{align*}
are isometric isomorphisms. Consequently, 
\begin{align*}
D^{(q)}(r|Y,\|\cdot\|_{[0,T]},p)& =
D^{(q)}(r|\pi_T^{(2)}(Y),\|\cdot\|_{[0,1]},p)\\
&\le c \,\IE[\|\pi^{(2)}_T (Y)\|_{\cH_1}^{\tilde p}]^{1/\tilde p}\
\frac 1{r+1}\\
&= c \,\IE[\|\sqrt T\, \pi^{(1)}_T (Y)\|_{\cH_1}^{\tilde p}]^{1/\tilde p}\
\frac 1{r+1}\\
&= c \,\sqrt T\, \IE[\| Y\|_{\cH_T}^{\tilde p}]^{1/\tilde p}\
\frac 1{r+1}.
\end{align*}
\end{proof}

\begin{lemma}
For any $p\ge1$, there exists a constant $c<\infty$ such that
$$
D^{(q)}(r|\bar A,\|\cdot\|_{[0,\tau]},p)\le c\,
(H^p(Z)+1)^{1/2p}\bigl[\sqrt T\frac 1r+ \IP(\tau>T)^{1/4p}\ \IE[\tau^{2p}]^{1/4p}\bigr]
$$
for all $T\ge0$, $r>0$, all $[0,\infty)$-valued r.v.'s $\tau$ and all side informations $Z$.
\end{lemma}

\begin{proof} Fix $T>0$ and $r\ge0$. The previous lemma and equation
  (\ref{eq2204-1}) imply
\begin{align*}
D^{(q,0)}(r|\bar A,\|\cdot\|_{[0,T]},p) &\le c_1 \ \sqrt T\  \IE[\|\bar A\|^{2p}]^{1/2 p}\ \frac 1{r+1}\\
&\le c_1\ \kappa_p^{1/2p}\ \sqrt T\ (H^{p}(Z)+1)^{1/2 p}\ \frac1{r+1}\\
&= {\textstyle \frac {c_2}2}\  \sqrt T\ (H^{p}(Z)+1)^{1/2 p}\ \frac1{r+1},
\end{align*}
for some appropriate constants $c_1, c_2$ depending only on $p$.
Consequently, there exists a codebook $\cC\subset\IC[0,\infty)$ of
size $\lfl e^r\rfl$ which contains
$0$ and satisfies
\begin{align}\label{eq2705-1}
\IE[\min_{\hat a\in\cC} \|\bar A-\hat  a \|_{[0,T]}^p]^{1/p}\le c_1 \ \sqrt T\, (H^{  p}(Z)+1)^{1/2 p}\ \frac1{r+1}.
\end{align}
Let $\hat A$ now denote an
$\|\cdot\|_{[0,\tau]}$-optimal
reconstruction for $\bar A$ in $\cC$. Then
\begin{align*}
\IE[\|\bar A-\hat A\|_{[0,\tau]}^p]^{1/p}&\le
\IE[1_{\{\tau\le T\}} \|\bar A-\hat A\|_{[0,\tau]}^p]^{1/p}+\IE[1_{\{\tau> T\}}
\|\bar A\|_{[0,\tau]}^p]^{1/p}\\
&=:I_1+I_2.
\end{align*}
It follows from equation (\ref{eq2705-1}) that
$$
I_1\le c_1 \ \sqrt T\, (H^{  p}(Z)+1)^{1/2 p}\, \frac1{r+1}.
$$
Moreover, the second term satisfies
\begin{align*}
I_2&\le \IE[ 1_{\{\tau> T\}} \sqrt\tau^p \|\bar
A\|_{\cH_\tau}^p]^{1/p}\\
&\le\IE[ 1_{\{\tau> T\}} \tau^p]^{1/2p}\,\IE[\|\bar A\|_{\cH}^{2p}]^{1/2p}\\
&\le \IP(\tau>T)^{1/4p} \,\IE[\tau^{2p}]^{1/4p}\,
\IE[\|\bar A\|_{\cH}^{2p}]^{1/2p}\\
&\le \kappa_{p}^{1/2p}\, \IP(\tau>T)^{1/4p}\, \IE[\tau^{2p}]^{1/4p}\, (H^p(Z)+1)^{1/2p}.
\end{align*}
Putting everything together yields the assertion.
\end{proof}

Let now $Z^{(r)}$ denote some side information depending on the rate
$r\ge0$ and let $W=\bar W^{(r)}+\bar A^{(r)}$ be the corresponding $(\cG_t)$-Doob-Meyer decomposition. Moreover, $\tau(r)$  $(r\ge0)$ denotes a random time in
$[0,\infty)$. We will need

\begin{lemma}\label{le1604-1} Let $p\ge1$ and $\alpha>0$. 
Suppose that for any $q<\infty$
$$
\sup_{r\ge0} \IE[\tau(r)^q]<\infty
$$
and 
\begin{align}\label{eq1604-3}
H^p(Z^{(r)})^{1/2p}=\cO(r^\alpha), \qquad r\to\infty.
\end{align}
Then for any $\eta>0$ and $\ep>0$
$$
D^{(q)}(r^\eta|\bar A^{(r)},\|\cdot\|_{[0,\tau(r)]},p)= \cO(r^{\alpha-\eta+\ep}),\qquad
r\to\infty.
$$
In particular, assumption (\ref{eq1604-3}) is fulfilled if the random variables
$Z^{(r)}$ have finite range and satisfy
\begin{align}\label{eq1604-4}
\log |\supp (Z^{(r)})|=\cO(r^{2\alpha}).
\end{align}
\end{lemma}

\begin{proof} The first statement is a consequence of the previous lemma.
In fact, for $\ep\in(0,\eta)$, $r\ge0$ and $T=T(r)=r^{2\ep}$ one has
$$
H^p(Z^{(r)})^{1/2p}\,\sqrt T\, \frac 1{r^\eta}= \cO(r^{\alpha-\eta+\ep}).
$$
Next, we set $q:=2p\frac{\eta-\ep}{\ep}$ and obtain
\begin{align*}
H^p(Z^{(r)})^{1/2p}\,& \IP(\tau(r) >T)^{1/4p}\,\IE[\tau(r)^{2p}]^{1/4p}\\
&\le
H^p(Z^{(r)})^{1/2p}\,
\Bigl(\frac{\IE[\tau(r)^q]}{T^q}\Bigr)^{1/4p}\,\IE[\tau(r)^{2p}]^{1/4p}\\
&= H^p(Z^{(r)})^{1/2p}\,\IE[\tau(r)^q]^{1/4p}
\,\IE[\tau(r)^{2p}]^{1/4p}\, r^{-(\eta-\ep)}\\
&=\cO(r^{\alpha-\eta+\ep}).
\end{align*}
It remains to show that (\ref{eq1604-4}) is a sufficient criterion for
(\ref{eq1604-3}).
Notice that for $p>1$, the map $x\mapsto (\log x)^p$ is not concave on the
  domain $[1,\infty)$. However, we can choose $c_1\ge0$ such that 
$$
f:[1,\infty)\to [0,\infty), x\mapsto (\log x)^p+c_1\,\log x
$$
is concave. Denote by $(q_x)$ the probability weights of some discrete r.v.\ $Z$. Then
\begin{align*}
H^p(Z)&=\sum_x q_x(\log(1/q_x))^p= \IE[(\log(1/q_Z))^p]\le
\IE[f(1/q_Z)]\\
&\le f\bigl(\IE[1/q_Z]\bigr)\le f(|\supp (Z)|).
\end{align*}
The assertion follows immediately.
\end{proof}

\section{Main representation of the diffusion}\label{sec7}

In this section we derive a representation of the diffusion which allows us to relate the coding complexity of the diffusion process
to that of the Wiener process. This will lead to the main theorems.

\begin{theo} \label{theo2705-1}
Fix $p\ge1$. There exist $\IC([0,\infty),\IR^d)$-valued r.v.'s $\bar X^{(r)}$, $\hat{\bar X}^{(r)}$ and $\bar W^{(r)}$, as well as a $\IC[0,\infty)$-valued r.v. $\hat\vphi^{(r)}$ $(r\ge0)$ such that 
\begin{itemize}
\item $\exists \gamma\in(0,1)$ s.th.\ $\log |\supp(\hat{\bar
    X}^{(r)},\hat\vphi^{(r)})|=\cO(r^{\gamma})$,
\item $\hat\vphi^{(r)}$ is a regular reconstruction for $\vphi$ with
  $\IE[\|\vphi-\hat\vphi^{(r)}\|_{[0,1]}]\longrightarrow0$
\item $X= \bar X^{(r)} +\bar W^{(r)}_{\hat \vphi^{(r)}(\cdot)}$,
\item $\bar W^{(r)}$ is a Wiener process that is independent of
  $\hat\vphi^{(r)}$
\item $\exists \delta>0$ s.th.\ $\IE[\|\bar X^{(r)}-\hat{\bar X}^{(r)}\|_{[0,1]}^p]^{1/p}=\cO(r^{-\frac12-\delta})$ as $r\to\infty$.
\end{itemize}
\end{theo}

\begin{proof}

1.) We let $\beta>0$ as in assumption (C) and set
$\alpha:=\beta/4$, $\delta_1:=\alpha/6$ and
$\gamma_1:=(1+\frac12\alpha)/(1+\alpha)<1$.
Due to Lemma \ref{le1604-2}, there exists some constant $c_1<\infty$
such that $\vphi$ admits $\IC[0,1]$-valued regular reconstructions $\hat\vphi^{(r)}$
$(r\ge0)$ satisfying
\begin{align}\label{eq1604-2}
|\supp (\hat\vphi^{(r)})|\le \exp\{r^{\gamma_1}\}\ \text{ and }\
\IE[\|\vphi-\hat\vphi^{(r)}\|_{[0,1]}^2]^{1/2}\le c_1\,\frac1{1+r^{1+3\delta_1}}
\end{align}
for all $r\ge0$. If $\hat\vphi^{(r)}$ is not monotonically increasing, we
replace $\hat\vphi^{(r)}$ by  the monotone function $t\mapsto \sup_{s\in[0,t]} \hat\vphi^{(r)}(s)$,
and condition (\ref{eq1604-2}) remains  valid.
By Lemma \ref{le1604-3}, it follows
\begin{align}\label{eq1027-3}
\IE[\|W_{\vphi(\cdot)}-W_{\hat\vphi(\cdot)}\|_{[0,1]}^p]^{1/p}=\cO(r^{-\frac12-\delta_1}).
\end{align}

2.) Let $W=\bar
W^{(r)}+\bar A^{(r)}$ be the Doob-Meyer decomposition of $W$ under the
side information $\hat \vphi^{(r)}$. 
We want to apply Lemma \ref{le1604-1}:
Let $\tau(r):=\hat\vphi^{(r)}(1)$. Then for any $\tilde p\in[1,\infty)$
\begin{align*}
\IE[\tau(r)^{\tilde p}]&=\IE[\bigl(\hat\vphi^{(r)}(1)\bigr)^{\tilde
  p}]\le \IE[\|\hat\vphi^{(r)}\|_{[0,1]}^{\tilde p}]\\
&\le \IE[\bigl(\|\hat\vphi^{(r)}-\vphi \|_{[0,1]}+\|\vphi\|_{[0,1]} \bigr)^{\tilde p}].
\end{align*}
Due to the regularity of $\hat\vphi^{(r)}$ it holds $ \|\vphi-
\hat\vphi^{(r)}\|_{[0,1]}\le \|\vphi\|_{[0,1]}$ and, hence,
$$
\IE[\tau(r)^{\tilde p}]^{1/\tilde p}\le 2\, \IE[\|\vphi\|_{[0,1]}^{\tilde
  p}]^{1/\tilde p}\le c_2<\infty
$$
is uniformly bounded for all $r\ge0$.
Recall that $|\supp (\hat\vphi^{(r)})|\le \exp\{r^{\gamma_1}\}$.
Moreover, for $\gamma_2:=(3+\gamma_1)/4<1$ it holds
$$
\gamma_1/2-\gamma_2=-\frac12-\frac{1-\gamma_1}4<-\frac12-\frac{1-\gamma_1}5=:-\frac12-\delta_2.
$$
Therefore Lemma \ref{le1604-1} implies that
$$
D^{(q)}(r^{\gamma_2}|\bar A^{(r)}, \|\cdot\|_{[0,\tau(r)]},
p)=\cO(r^{-\frac12-\delta_2}),
$$
and there exist a constant $c_3<\infty$
and continuous reconstructions $\hat{\bar A}^{(r)}$ such that
\begin{align}\label{eq1027-2}
\IE[\|\bar A^{(r)}-\hat{\bar A}^{(r)}\|_{[0,\tau(r)]}^p]^{1/p}\le c_3
r^{-\frac12-\delta_2}\ \text{ and }\  \log|\supp (\hat{\bar A}^{(r)})|\le r^{\gamma_2}.
\end{align}
Finally notice that by Lemma \ref{le0705-1}, there exists a constant
$c_4<\infty$ and reconstructions $\hat A^{(r)}$ for $A$ such that for all $r\ge0$
\begin{align}\label{eq0705-1}
\IE[\|A-\hat A^{(r)}\|^p_{[0,1]}]^{1/p} \le c_4
\frac1{r^{\frac12+\delta_3}}\ \text{ and }\  \log|\supp (\hat A^{(r)})|\le
r^{\gamma_3},
\end{align}
where $\delta_3:=\frac16$ and $\gamma_3=\frac23$.\\
3.) We rewrite $X$ in terms of the new r.v.'s:
\begin{align*}
X_t&=A_t+M_t=A_t+W_{\vphi(t)}=A_t+(W_{\vphi(t)}-W_{\hat\vphi^{(r)}(t)})+W_{\hat\vphi^{(r)}(t)}\\
&=\underbrace{A_t+(W_{\vphi(t)}-W_{\hat\vphi^{(r)}(t)})+\bar
  A^{(r)}_{\hat\vphi^{(r)}(t)}}_{=:\bar X^{(r)}_t}+\bar W^{(r)}_{\hat\vphi^{(r)}(t)}.
\end{align*}
Due to (\ref{eq1027-3}), (\ref{eq1027-2}) and (\ref{eq0705-1}) it follows that the process $\hat{\bar X}^{(r)}_t:=\hat A_t^{(r)}+\hat{\bar
  A}^{(r)}_{\hat\vphi^{(r)}(t)}$ satisfies for
  $\delta:=\min(\delta_1,\delta_2,\delta_3)>0$ and $\gamma=\max(\gamma_1,\gamma_2,\gamma_3)<1$,
\begin{align*}
\IE[\|\bar X^{(r)}- \hat{\bar X}^{(r)}\|_{[0,1]}^p]^{1/p}&\le \IE[\|A-\hat
A^{(r)}\|_{[0,1]}^p]^{1/p} +
\IE[\|W_{\vphi(\cdot)}-W_{\hat\vphi(^{(r)}\cdot)}\|_{[0,1]}^p]^{1/p}\\
& \qquad +\IE[\|\bar
A^{(r)}_{\hat\vphi^{(r)}(\cdot)}-\hat{\bar A}^{(r)}_{\hat\vphi^{(r)}(\cdot)}\|_{[0,1]}^p]^{1/p}=\cO(r^{-\frac12-\delta})
\end{align*}
and
$$
\log | \supp (\hat{\bar X}^{(r)},\hat\vphi^{(r)})|\le
 r^{\gamma_1}+r^{\gamma_2}+r^{\gamma_3}= \cO(r^{\gamma}).
$$
\end{proof}

\section{The quantization complexity of $X$ in $\IC([0,1],\IR^d)$}

We are now in a position to prove Theorem \ref{theo1109-1}.
Fix $p\ge1$ arbitrary and let  $K\in[\frac{\pi}{\sqrt 8},\pi]$ such that
\begin{align}\label{eq3006-2}
D^{(q)}(r|W,\|\cdot\|_{[0,1]},p)\sim K \,\frac 1{\sqrt r},\qquad r\to\infty.
\end{align}

\begin{lemma}
\begin{align*}
D^{(q)}(r|X,\|\cdot\|_{[0,1]},p) \lesssim K\, \IE[
\|\sig\|_{L^2[0,1]}^p]^{1/p}\  \frac {1}{\sqrt r }.
\end{align*}
\end{lemma}

\begin{proof}Let $\bar X=\bar X^{(r)}$, $\hat {\bar X}=\hat {\bar X}^{(r)}$, $\hat\vphi=\hat\vphi^{(r)}$ and
$\bar W=\bar W^{(r)}$  be as in Theorem \ref{theo2705-1}, and let $\tau=\hat\vphi(1)$.
There exists a discrete r.v.\ $\hat W$ such that,
conditional upon $\hat\vphi$, the random law $\cL(\hat W|\hat
\vphi)$ is supported on a set of size $\lfl e^r\rfl$ and 
$$
\IE[\|\bar W-\hat W\|_{[0,\tau]}^p|\hat\vphi] = \sqrt\tau ^p D^{(q)}(r|W,\|\cdot\|_{[0,1]},p)^p.
$$
Then
\begin{align*}
\IE[\|\bar W-\hat W\|_{[0,\tau]}^p]^{1/p}= \IE[ \sqrt\tau ^p ]^{1/p}\, D^{(q)}(r|W,\|\cdot\|_{[0,1]},p).
\end{align*}
Notice that  $\lim_{r\to\infty} \IE[\sqrt\tau
^p]^{1/p}=\IE[\|\sig\|_{L^2[0,1]}^p]^{1/p}$. As a reconstruction, we choose $\hat X:=\hat X^{(r)}:= \hat{\bar X}+\hat W$,
so that
\begin{align*}
\IE[\|X-\hat X\|_{[0,1]}^p]^{1/p}&\le \IE[\|\bar X-\hat{\bar
  X}\|_{[0,1]}^p]^{1/p}+ \IE[\|\bar W-\hat W\|_{[0,\tau]}^p]^{1/p}\\
&\lesssim \IE[\|\sig\|_{L^2[0,1]}^p]^{1/p}\, K\,\frac1{\sqrt r}.
\end{align*}
Moreover, $\hat X$ has range of size $e^{(1+o(1))r}$.
\end{proof}

Now we turn to the proof of the converse inequality:
\begin{lemma}
\begin{align*}
D^{(q)}(r|X,\|\cdot\|_{[0,1]},p) \gtrsim K\, \IE[
\|\sig\|_{L^2[0,1]}^p]^{1/p}\  \frac {1}{\sqrt r }.
\end{align*}
\end{lemma}

\begin{proof} Let $\bar X$, $\hat {\bar X}$, $\hat\vphi$ and $\bar W$
  be as in Theorem \ref{theo2705-1}. 
Denote by $\hat X=\hat X^{(r)}$ an arbitrary reconstruction for $X$
that has range of size $\lfl e^r\rfl$.
Let 
$$
\hat W_t:=\hat W^{(r)}_t:=\hat X_{\hat\vphi^{-1}(t)}- \hat{\bar
  X}_{\hat\vphi^{-1}(t)}\qquad (t\in[0,\tau]),
$$ 
where $\hat \vphi^{-1}(t):=\inf\{s\ge0: \hat\vphi(s)\ge t\}$. 
Since $\bar W_t=X_{\hat\vphi^{-1}(t)}-\bar X_{\hat\vphi^{-1}(t)}$ for $t\in[0,\tau]$ one has
\begin{align}\begin{split} \label{eq3006-1}
\IE[\|\bar W-\hat
W\|_{[0,\tau]}^p]^{1/p}&\le \IE[\|X_{\hat\vphi^{-1}(\cdot)}-\hat
X_{\hat\vphi^{-1}(\cdot)}\|_{[0,\tau]}^p]^{1/p}+  \IE[\|\bar
X_{\hat\vphi^{-1}(\cdot)}-\hat {\bar
  X}_{\hat\vphi^{-1}(\cdot)}\|_{[0,\tau]}^p]^{1/p}\\
&\le  \IE[\|X-\hat
X\|_{[0,1]}^p]^{1/p}+\IE[\|\bar X-\hat{\bar X}\|_{[0,1]}^p]^{1/p}.
\end{split}
\end{align}
On the other hand,  the random law $\cL(\hat W|\hat\vphi)$ is
supported on a set of size $|\supp (\hat {\bar X})|\cdot \lfl e^r
\rfl= \exp\{r+\cO(r^\gamma)\}$ for some $\gamma\in(0,1)$,
 and given $\hat\vphi$ the process $\bar W$ is a Wiener process. Therefore,
$\IE[\|\bar W-\hat W\|_{[0,\tau]}^p|\hat\vphi ]\ge
K^p\,{\sqrt\tau}^p\,(1-o(1))/{\sqrt r}^p$ as $r\to\infty$. Here, the
$o(1)$-term depends only on $r$ but not on the realization of $\hat\vphi$.
Consequently,
\begin{align*}
\IE[\|\bar W-\hat
W\|_{[0,\tau]}^p]^{1/p}&=\IE\bigl[\IE[\|\bar W-\hat
W\|_{[0,\tau]}^p|\hat\vphi ]\bigr]^{1/p}\\
&\gtrsim \IE\Bigl[K^p {\sqrt\tau}^p \frac1{\sqrt r ^p}
\Bigr]^{1/p}=K\, {\IE[{\sqrt\tau}^p]^{1/p}} \frac 1{\sqrt r}.
\end{align*}
Note that $\lim_{r\to\infty} \IE[\sqrt
{\tau}^p]^{1/p}=\IE[\|\sig\|_{L^2[0,1]}^p]^{1/p}$ and, hence, 
(\ref{eq3006-1}) implies that
\begin{align*}
\IE[\|X-\hat
X\|_{[0,1]}^p]^{1/p}&\ge (1-o(1))\,K\,\IE[\|\sig\|_{L^2[0,1]}^p]^{1/p}
\frac 1{\sqrt r}-\IE[\|\bar X-\hat{\bar X}\|_{[0,1]}^p]^{1/p}\\
&\sim K\,\IE[\|\sig\|_{L^2[0,1]}^p]^{1/p}
\frac 1{\sqrt r}.
\end{align*}
\end{proof}

\section{Coding complexity of $W$ in general $L^p(\nu)$-spaces}

Let $p,q\ge1$, $T>0$ and $\nu\in\cM[0,T)$, where $\cM[0,T)$ denotes the set of finite positive measures on the Borel sets of $[0,T)$. The objective of this section is to provide an upper estimate for $D^{(q)}(r|W,\|\cdot\|_{L^q(\nu)},p)$.

\begin{lemma} \label{le1504-1} Let $p\ge q\ge1$. There exists a constant
  $c=c(p)$ such that for all  $r, \Delta r, T \ge0$ and
  $\nu\in \cM([0,T))$
\begin{align*}
D^{(q)}(r+&\Delta r|W,\|\cdot\|_{L^q(\nu)},p)\\
&\le  \nu[0,T)^{1/q} \sqrt T \bigl[ D^{(q)}(r|W,\|\cdot\|_{L^q[0,1)},p)+c\, e^{-\Delta r/2d}\bigr].
\end{align*}
\end{lemma}

The proof of the lemma uses a simple consequence of \cite{Der03}, Proposition 6.3.4 (see also Corollary 6.7 in \cite{GraLu00}):

\begin{propo}\label{prop0730-1}
Let $\tilde p> p\ge 1$ and $m\iN$. There exists a constant $c<\infty$ such that for all $\IR^m$-valued random variables $Z$ and all $r\ge0$, 
$$
D^{(q)}(r|Z,|\cdot|,p)\le c\, \IE[|Z|^{\tilde p}]^{1/\tilde p}\, e^{-r/m}.
$$
\end{propo}

\begin{proof}[ of Lemma \ref{le1504-1}] First consider the case where $\nu$ is a probability measure on $[0,T)$.
For $t\in[0,T)$ let 
$$
\th_t:[0,T)\to[0,T):=s+t\text{ mod }T 
$$
and let $\nu_t:=\nu\circ \th_t^{-1}$ denote the measure induced by the map
$\th_t$. Now let $\hat W$ be a r.v.\ attaining at most $e^r$ different values with
$$
\IE[\|W-\hat W\|_{L^q[0,T)}^p]^{1/p}=D^{(q)}(r|W,\|\cdot\|_{L^q[0,T)},p).
$$
Using the
Jensen inequality we arrive at
\begin{align*}
&\int_0^T\IE\Bigl[\Bigl(\int_{[0,T)} |W_u-\hat W_u|^q\, d\nu_t(u)
\Bigr)^{p/q}\Bigr]\,\frac{dt}{T}\\
& \ \ \le  \IE\Bigl[ \Bigl(\int_0^T \int_{[0,T)} |W_u-\hat W_u|^q\, d\nu_t(u)
\,\frac{dt}{T}\Bigr)^{p/q}\Bigr]\\
& \ \ = \IE\Bigl[ \Bigl(\int_0^T |W_u-\hat W_u|^q\, \frac{du}T
\Bigr)^{p/q}\Bigr]=\frac 1{T^{p/q}}\,D^{(q)}(r|W,\|\cdot\|_{L^q[0,T]},p)^p
\end{align*}
In particular, there exists $t\in[0,T)$ such that
\begin{align*}
\IE[\|W_{\th_t(\cdot)}-\hat W_{\th_t(\cdot)}\|_{L^q(\nu)}^p]^{1/p}&=\IE\Bigl[\Bigl(\int_{[0,T)} |W_u-\hat W_u|^q\, d\nu_t(u)
\Bigr)^{p/q}\Bigr]^{1/q}\\
& \le \frac 1{T^{1/q}}\,D^{(q)}(r|W,\|\cdot\|_{L^q[0,T)},p).
\end{align*}
Using that
$D^{(q)}(r|W,\|\cdot\|_{L^q[0,T)},p)=T^{\frac12+\frac1q}\,D^{(q)}(r|W,\|\cdot\|_{L^q[0,1)},p)$
one obtains
$$
\IE[\|W_{\th_t(\cdot)}-\hat W_{\th_t(\cdot)}\|_{L^q(\nu)}^p]^{1/p}\le
\sqrt{T}\,D^{(q)}(r|W,\|\cdot\|_{L^q[0,1)},p).
$$

Now let for $s\in[0,T)$
\begin{align*}
W_s':=\begin{cases}
W_{\th_t(s)}-W_t & \text{ if }s+t<T\\
W_{\th_t(s)}+W_1-W_t & \text{ if }s+t\ge T.
\end{cases}
\end{align*}
Clearly, $(W'_s)_{s\in[0,T)}$ is a standard $d$-dimensional Wiener process. We
consider reconstructions $(\hat W'_s)$ of the form
\begin{align}\label{eq0704-1}
\hat W_s':=\begin{cases}
\hat W_{\th_t(s)}+A_1 & \text{ if }s+t<T\\
\hat W_{\th_t(s)}+A_2 & \text{ if }s+t\ge T,
\end{cases}
\end{align}
where $A_1$ and $A_2$ denote $\IR^d$-valued r.v.'s which still need to
be fixed appropriately.
Then the coding error can be controlled by
\begin{align}\label{eq0704-2}
\|W'-\hat W'\|_{L^q(\nu)}\le \|W_{\th_t(\cdot)}-\hat W_{\th_t(\cdot)}\|_{L^q(\nu)}
+\Bigl|\begin{pmatrix} -W_t-A_1  \\ W_1-W_t-A_2 \end{pmatrix}\Bigr|.
\end{align}
Set $Z:=(-W_t,W_1-W_t)^t$. Clearly, $\IE[|Z|^{2p}]$ is finite and due
to Proposition \ref{prop0730-1} there exists a constant
$c_1=c_1(p)<\infty$ only depending on $p\ge1$ 
such that for all $\Delta r\ge0$
$$
D^{(q)}(\Delta r|Z,|\cdot|,p) \le c_1\, \IE[|Z|^{2p}]^{1/2p}
\,e^{-\Delta r/2d}.
$$
Notice that $\IE[|Z|^{2p}]^{1/2p}\le  c_2\,\sqrt T$ for some
constant $c_2=c_2(p)<\infty$ only depending on $p$.
Choose now $A_1,A_2$ such that the random variable $(A_1,A_2)$ has range of size $e^{\Delta r}$ and satisfies 
$$
\IE\Bigl[\Bigl| \begin{pmatrix}-W_t-A_1\\ W_1-W_t-A_2
\end{pmatrix}\Bigr|^p\Bigr]^{1/p}\le c \, \sqrt T\,e^{-\Delta r/2d},
$$
where $c:=2\,c_1\,c_2$.
 Then $\hat W'$ as defined
in (\ref{eq0704-1}) has range of size $e^{r+\Delta r}$ and by (\ref{eq0704-2}) it holds
$$
\IE[\|W'-\hat W'\|^p_{L^q(\nu)}]^{1/p}\le \sqrt{T}
D^{(q)}(r|W,\|\cdot\|_{L^q[0,1)},p)+c\,\sqrt T\, e^{-\Delta r/2d}.
$$ 
For a general finite measure $\nu\in\cM[0,T)$ there exists a probability measure $\nu_1\in\cM[0,T)$ such that $\nu=\nu[0,T) \cdot \nu_1$, and one has
$$
\|W'-\hat W'\|_{L^q(\nu)}=\nu[0,T)^{1/q}\, \|W'-\hat W'\|_{L^q(\nu_1)}.
$$
The assertion follows immediately.
\end{proof}


\section{The quantization complexity of $X$ in $L^p[0,1]$}

For fixed $p\in[0,\infty)$, we consider asymptotic
upper bounds for 
$$
D^{(q)}(r|X,\|\cdot\|_{L^p[0,1]},p).
$$
We denote by $K=K(p)<\infty$ a constant with the property that
\begin{align}\label{eq2204-2}
D^{(q)}(r|W,\|\cdot\|_{L^p[0,1]},p)\lesssim K\frac 1{\sqrt r},\qquad r\to\infty.
\end{align}
Due to \cite{DFMS03} such a constant exists.

\begin{theo}\label{theo0306-1}
\begin{align}\label{eq1604-1}
D^{(q)}(r|X,\|\cdot\|_{L^p[0,1]},p) \lesssim  K\,\IE[ \|\sig\|_{L^{2p/(p+2)}[0,1]}^p]^{1/p}\, \frac {1}{\sqrt r }.
\end{align}
\end{theo}

\begin{proof}
Let $\bar X=\bar X^{(r)}, \hat{\bar X}=\hat{\bar X}^{(r)},\hat \vphi= \hat \vphi^{(r)}$ and $\bar
W=\bar W^{(r)}$ be as in Theorem \ref{theo2705-1}.
Recall that  the time change $\hat\vphi=\hat\vphi^{(r)}$ and the Wiener process $\bar
W=\bar W^{(r)}$ are independent.  We fix $n\iN$ and denote 
$\tau_i:=\tau_i(r):=\hat\vphi^{(r)}(i/n)$ ($i=0,\dots,n$) and $\tau:=\tau(r):=\hat\vphi^{(r)}(1)$.
Moreover, let for $t\in[\tau_{i-1},\tau_i)$,
$$
W'_{t} := \bar W_t- \bar W_{\tau_{i-1}}\ \ \text{ and }\ \ W''=\bar W_{\tau_{i-1}}.
$$
Clearly, $\bar W=W'+W''$.
We define reconstructions for $W'$ and $W''$ separately.

In the following, let $\Delta r=\sqrt r$.
Set $Y:=(W_{\tau_1},\dots,W_{\tau_n})^t$.
Assume first that $(\tau_i)$ is deterministic. Then one has
$$
D^{(q)}(\Delta r|W'',\|\cdot\|_{[0,\tau]},p)\le D^{(q)}(\Delta r|Y,|\cdot|,p).
$$
As in the proof of Lemma \ref{le1504-1} we conclude that there exists a universal constant $c_1=c_1(p)<\infty$ such that
$$
 D^{(q)}(\Delta r|Y,|\cdot|,p)\le \textstyle{\frac {c_1}2}\, \sqrt{\tau}\, e^{-\Delta r/nd}.
$$
Now let $(\tau_i)$ be random as before. Since $W''$ is independent of $\hat\vphi$ there exists a discrete
r.v.\ $\hat W''$ such that conditional upon $\hat\vphi$,
$\hat W''$ has range of size $\lfl e^{\Delta r}\rfl$ and 
it holds
$$
\IE[ \|W''-\hat W''\|_{[0,\tau]}^p|\hat\vphi]^{1/p} \le c_1\, \sqrt{\tau}\, e^{-\Delta r/nd}.
$$
Consequently,
$$
\IE[ \|W''-\hat W''\|_{[0,\tau]}^p]^{1/p} \le c_1\, \IE[ \sqrt\tau ^p]^{1/p}\, e^{-\Delta r/nd}.
$$
Clearly, this expression is of order $o(1/\sqrt r)$ since
$\IE[\sqrt\tau ^p]$ is uniformly bounded.

Now we construct a reconstruction for $W'$.
We decompose the random measure $\nu:=\nu^{(r)}:=\lam (\hat\vphi^{-1})$ into the sum 
$$
\nu=\sum_{i=0}^n
\nu_i,
$$
 where $\nu_i=\nu|_{(\tau_{i-1},\tau_{i+1})}$ for $i=1,\dots, n$, and
 $\nu_0$ contains the remaining mass of $\nu$.
Then $\nu_i[\tau_{i-1},\tau_i)\le 1/n$.
We set $\Delta \tau_i=\Delta \tau_i^{(r)}=\tau_i-\tau_{i-1}$ for $i=1,\dots,n$, and choose
$$
r_i=r_i(r,\hat\vphi)=\frac{\Delta\tau_i^{p/(p+2)}}{\sum_{j=1}^n \Delta\tau_j^{p/(p+2)}}\, r\vee
\sqrt r.
$$
Due to Lemma \ref{le1504-1}, 
there exists a constant $c_2=c_2(p)$ and  reconstructions $\hat W'=\hat W'^{(r)}$ such that conditional upon
$\hat \vphi$, $\hat W'$ has range of size
$\exp\{r_1+\dots+r_n+2dn\Delta r\}$
and
\begin{align*}
\IE[\| W'-\hat W'\|_{L^p(\nu)}^p|\hat\vphi]&=\sum_{i=1}^n \IE[\|W'-\hat W'\|^p_{L^p(\nu_i)}|\hat\vphi]\\
&\le \frac 1n\sum_{i=1}^n
\Delta \tau_i^{p/2}\bigl[D^{(q)}(r_i|W,\|\cdot\|_{L^p[0,1)},p)+c_2 e^{-\Delta r}\bigr]^p.
\end{align*}

Notice that $r_i$ converges uniformly to $\infty$ (independently of
the realization of $\hat\vphi$). Thus, using condition (\ref{eq2204-2}), we conclude that
\begin{align*}
\Delta\tau_i^{p/2} \,D^{(q)}(r_i|W,\|\cdot\|_{L^p[0,1)},p)^p \le 
(K+o(1))^p \, \Delta\tau_i^{p/2} \frac 1{\sqrt {r_i}^p}\\
\le (K+o(1))^p\, \Delta\tau_i^{p/(p+2)} \bigl(\sum_{j=1}^n \Delta
\tau_j^{p/(p+2)}\bigr)^{p/2}\frac1{\sqrt r ^p}.
\end{align*}
Here, the $o(1)$-term converges uniformly to $0$ as $r\to\infty$, independently of the
realization of $\hat\vphi$. Since
$\IE[\tau_i^{p/2}]\le\IE[\tau^{p/2}]$ is uniformly bounded, it follows that
$$
\IE[\| W'-\hat W'\|_{L^p(\nu)}^p]\le \frac {(K+o(1))^p}{\sqrt r ^p}\IE
\Bigl[\frac1n\bigl(\sum_{j=1}^n \Delta
\tau_j^{p/(p+2)}\bigr)^{(p+2)/2} \Bigr].
$$
Note that for all $i=1,\dots,n$, $\Delta\tau_i$ converges in probability to $\vphi(i/n)-\vphi((i-1)/n)$ and, hence, by dominated convergence it follows that 
$$
\IE[\| W'-\hat W'\|_{L^p(\nu)}^p]^{1/p}\le (K+o(1)) \,Z_n \,\frac1{\sqrt r },
$$ 
where
$$
Z_n:=\IE\Bigl[\frac1n\bigl(\sum_{j=1}^n  (\vphi(i/n)-\vphi((i-1)/n))^{p/(p+2)}\bigr)^{(p+2)/2} \Bigr]^{1/p}.
$$
Therefore, the reconstruction $\hat W:=\hat W^{(r)}:=\hat W'+\hat W''$ satisfies 
$$
\IE[\|\bar W-\hat W\|_{L^p(\nu)}^p]^{1/p} \lesssim  K \, Z_{n} \,
\frac1{\sqrt r}.
$$

It remains to combine all estimates to control the quality and
complexity of the
reconstruction
$$
\hat X:=\hat X^{(r)}:=\hat{\bar X}+\hat W_{\hat \vphi(\cdot)}.
$$
One has:
\begin{align*}
\log | \supp (\hat X)| &\le
\log | \supp (\hat {\bar X},\hat W,\hat
\vphi)|\\
 & \le
r+(2dn+1)\Delta r+ \log | \supp (\hat {\bar X},\hat
\vphi)|\\
&= (1+o(1))\ r.
\end{align*}
Moreover, notice that
\begin{align*}
\|X-\hat X\|_{L^p[0,1]}&\le \|\bar X-\hat {\bar
  X}\|_{[0,1]}+ \|\bar W_{\hat \vphi(\cdot)}-\hat  W_{\hat \vphi(\cdot)}\|_{L^p[0,1]}\\
&= \|\bar X-\hat {\bar
  X}\|_{[0,1]}+\|\bar W-\hat  W\|_{L^p(\nu)},
\end{align*}
hence:
$$
\IE[\|X-\hat X\|_{L^p[0,1]}^p]^{1/p}\lesssim K\ Z_n\ \frac 1{\sqrt r}.
$$
The statement is valid for all $n\iN$ and it remains to show that $\lim_{n\to\infty} Z_n=
\IE[\|\sig\|_{L^{2p/(p+2)}[0,1]}^p]^{1/p}$.
Let for fixed $n\iN$ and for $i=1,\dots,n$ and $t\in[(i-1)/n,i/n)$, 
$$
\bar\sig^2_t= n \int_{(i-1)/n}^{i/n} \sig_t^2\,dt.
$$
Then we can rewrite $Z_n$ in terms of $\bar\sig^2$:
$$
Z_{n}=\IE \Bigl[\Bigl(\int_0^1  (\bar\sig^2_t)^{p/(p+2)}\, dt\Bigr)^{(p+2)/2} \Bigr]^{1/p}.
$$
As $n$ tends to infinity, $\bar\sig^2$ converges pointwise  to $\sig^2$. Hence, the result follows by the dominated convergence theorem.
\end{proof}

\bibliography{biblio}
\bibliographystyle{ims}

\end{document}